\documentclass{amsart}
\usepackage{amssymb}

\newtheorem{theorem}{Theorem}[section]

\newtheorem{lemma}[theorem]{Lemma}
\newtheorem{proposition}[theorem]{Proposition}

\newtheorem{remark}{Remark}[section]

\newtheorem{definition}{Definition}[section]

\theoremstyle{definition}
\theoremstyle{remark}
\numberwithin{equation}{section}



\begin{document}

\title{Operators on $C(\omega ^\alpha )$ which do not preserve
$C(\omega^\alpha )$}
\author{Dale E. Alspach}
\address{Oklahoma State University\\
Department of Mathematics\\
Stillwater, OK 74078
}
\date{October 14, 1996}
\email{alspach@math.okstate.edu}
\begin{abstract}
It is shown that if $\alpha ,\zeta $ are ordinals such that $1\leq \zeta
<\alpha <\zeta \omega ,$ then there is an operator from $C(\omega ^{\omega
^\alpha })$ onto itself such that if $Y$ is a subspace of $C(\omega ^{\omega
^\alpha })$ which is isomorphic to $C(\omega ^{\omega ^\alpha })$ $,$ then
the operator is not an isomorphism on $Y.$ This contrasts with a result of
J. Bourgain that implies that there are uncountably many ordinals $\alpha $
for which any operator from $C(\omega ^{\omega ^\alpha })$ onto itself there
is a subspace of $C(\omega ^{\omega ^\alpha })$ which is isomorphic to $%
C(\omega ^{\omega ^\alpha })$ on which the operator is an isomorphism.
\end{abstract}

\maketitle

In an earlier paper \cite{A1} we proved that there is an operator on $%
C(\omega ^{\omega ^2})$ which is not an isomorphism on any subspace which is
isomorphic to $C(\omega ^{\omega ^2})$ but the operator is onto $C(\omega
^{\omega ^2}).$ This in contrast with the situation for $C(\omega )$ and $%
C(\omega ^\omega )$ where there are no surjective operators which do not
preserve isomorphically a copy of the space,
\cite{P}, \cite{A2}. Subsequently, Bourgain 
\cite{B} proved a very general result which gives an estimate on the size of
the ordinal $\beta $ such that any operator on $C(\omega ^{\omega ^\alpha })$
which is surjective must be an isomorphism on a subspace isomorphic to $%
C(\omega ^{\omega ^\beta })$. Recently Gasparis \cite{Gasp} has generalized
the example in \cite{A1} to the case of operators on $C(\omega ^{\omega
^{\alpha +1}})$ to show that there are surjective operators on these spaces
which do not preserve a copy of $C(\omega ^{\omega ^{\alpha +1}}).$ For most
ordinals $\alpha $ this is very far from the estimate given by Bourgain. 

Bourgain used the Szlenk index and a combinatorial argument in the proof of
his result. Implicit in his proof is the notion of $\gamma $-families of
sets which was independently developed by Wolfe, \cite{wolfe}, and in \cite
{A1}. The existence of $\gamma$-families with asssociated measures is an
indication of the amount of topological disjointness in a subset of $C(K)^*$
 whereas the Szlenk
index only indicates disjointness. Bourgain essentially shows
that a large Szlenk index forces the
the existence of a $\gamma$-family of sets with the size of $\gamma$
dependent on the Szlenk index.
The existence of a $\gamma$-family is equivalent to a condition 
on an ordinal index which we have
named the Wolfe index. Thus from this view point
Bourgain proves that the Szlenk index gives
some lower bound on the Wolfe index. In some cases he obtains that
the two indices are of
roughly equivalent size. 
In this paper we give a very general construction of
examples of the type in \cite{A1} and \cite{Gasp} and show that there are
many more ordinals for which the Szlenk and Wolfe index are very different.

We will use notation similar to that 
in \cite{A1}. In particular, if $\gamma $ is
an ordinal, $C(\gamma )$ is the space of continuous functions on the
ordinals less than or equal to $\gamma $ with the order topology, which we
denote by $[1,\gamma ].$ If $K$ is a topological space, $K^{(\beta )}$ is
the $\beta $-derived set of $K.$ If $L\subset C(K)^*$, then $L^{(\beta)}$,
will be the $\beta$ derived set of $L$ 
with respect to the $\text{w}^*$-topology.
 If $K$ is a countable compact Hausdorff
space, then $K$ is homeomorphic to $[1,\omega ^\beta n],$ where the
cardinality of $K^{(\beta )}$ is $n,$ for some $n\in \mathbb{N},$ \cite{MS}.
It was shown in \cite{BP} that $C(\omega ^{\omega ^\alpha })$ is isomorphic
to $C(\omega ^\beta n)$ if and only if $\omega ^\alpha \leq \beta <\omega
^{\alpha +1}.$ Thus from the point of view of the isomorphic theory of
Banach spaces, the spaces $C(\omega^{\omega^\alpha})$, $\alpha<\omega_1$,
are a complete set of representatives of the $C(K)$-spaces for $C(K)$
separable and $K$ countable.

\section{A Topological Construction}

In order to define the operator we need to develop a method of constructing
special sets of measures on $\omega ^{\omega ^\alpha }$ which are
homeomorphic to $\omega ^{\omega ^\alpha }$ but which have supports which
are almost
disjoint but are not topologically well separated. In \cite{A1} we used
the porcupine topology, \cite{BdL}, to effect the construction.
Here we use a similar
construction but with somewhat different notation. The operators that we
construct are of the same form as that in our earlier work. Namely, we will
produce 
a compact Hausdorff space $K$ and a $w^{*}$-closed subset $L$ of $C(K)^{*}$
and we define an operator from $C(K)$ into $C(L)$ by evaluation. In this
paper we need to iterate the construction in \cite{A1}. To this end we
introduce a general procedure for extending a pair $(K,L)$ by a sequence of
spaces $K_n$, where $K$ and $K_n$ are compact Hausdorff spaces, each $K_n$
has a distinguished point $k_{n,0}$ and $L$ is a set of purely atomic
finitely supported probability measures on $K$.

 For each $k\in K$, we let $%
L(k)=\{l\in L|l(k)\neq 0\}\cup \{\emptyset \}$ and $S(k,L)$ be the one point
compactification of $\sum_n\sum_{l\in L(k)}K_n\setminus \{k_{n,0}\}$, where
we use $\sum_{i\in I}W_i$ to denote the disjoint sum of topological spaces $%
W_i$ with the topology generated by sets of the form $\cup _{i\in I}G_i$
with $G_i$ open for each $i\in I$. We denote the points of $S(k,L)$ as $4$%
-tuples $(k,l,n,j)$ where $l\in L(k)$ and $j\in K_n$ .The point added will
be denoted $(k,\emptyset )$ although it is also $(k,l,n,k_{n,0})$ for any $%
l\in L(k)$ and $n$. Note that if $L(k)=\{\emptyset \},$ 
$S(k,L(k))=\{(k,\emptyset )\}.$ We want to define a topology on the disjoint
union of the sets $S(k,L)$. Intuitively we want to glue $S(k,L)$ to $K$ at
the point $k$ by identifying $(k,\emptyset )$ with $k$. We also want to
extend the measures $L$ by sets of measures $L_n$ on $K_n$ and make a copy
of $L_n$ for each $l\in L(k),n\in \mathbb{N}$. More formally we make the
following definition.

\begin{definition}
Suppose that $K$ and $K_n,n\in \mathbb{N,}$ are compact Hausdorff spaces, $K_n$
has a distinguished point $k_{n,0}$ and $L,L_n,$ are sets of purely atomic
disjointly supported probability measures on $K$, $K_n$,respectively, with $%
\delta _{k_{n,0}}\in L_n$, for each $n.$ Define $(K,L)\bigotimes
\{(K_n,L_n),n\in \mathbb{N}\}$ to be the pair $(K^{\prime },L^{\prime })$ where 
$K^{\prime }$ is the compact Hausdorff space and $L^{\prime }$ is the set of
atomic probability measures on $K^{\prime }$ described below. $K^{\prime }$
is the set of 4-tuples $(k,l,n,j_n)$, $k\in K$, $l\in L(k),n\in \mathbb{N}$ and 
$j_n\in K_n$ with the topology generated by sets of the form 
\begin{eqnarray*}
\cup _{k\in K}G_k\cup \cup _{k\in G}\{(k,l,n,j_n) &:&k\in K,l\in L(k),n\in 
\mathbb{N,}j_n\in K_n\} \\
&&\setminus \cup _{(k,l,n)\in F}\{k\}\times \{l\}\times \{n\}\times
F_{k,l,n},
\end{eqnarray*}
where $G_k$ is an open subset of 
\[
\{(k,l,n,j_n):l\in L(k),n\in \mathbb{N,}j_n\in K_n\setminus
\{k_{n,0}\}\}=S(k,L)\setminus \{(k,\emptyset )\}
\]
for each $k$, $G$ is an open set in $K$, $F$ is a finite set of triples $%
(k,l,n)$ with $k\in K,n\in \mathbb{N}$ and $l\in L(k)$, and $F_{k,l,n}$ is a
compact subset of $K_n\setminus \{k_{n,0}\}.$ For each $k\in K$ we identify
all of the points $(k,l,n,k_{n,0})$ such that $l\in L(k),n\in \mathbb{N}$ with
the point $(k,\emptyset )$. (Formally $K^{\prime }$ is a set of equivalence
classes of 4-tuples, but only the elements with fourth entry $k_{n,0}$ are
in non-trivial classes.) Let $\phi $ be the map from $K$ into $K^{\prime }$
defined by $\phi (k)=(k,\emptyset )$ and let $\Phi $ be the map from $%
\mathcal{M}(K)$ into $\mathcal{M}(K^{\prime })$ which is induced by $\phi $.
Let 
\[
L^{\prime }=\sum_{k\in \text{supp }l}l(k)\sum_{j_n\in K_n}l_n(j_n)\delta
_{(k,l,n,j_n)}):l\in L,n\in \mathbb{N},l_n\in L_n\}.
\]
\end{definition}

In keeping with our identification 
\[
\sum_{k\in \text{supp }l}l(k)\delta _{(k,l,n,k_{n,o})}=\sum_{k\in \text{supp 
}l}l(k)\delta _{(k,\emptyset )}=\Phi (l) 
\]
for each $l\in L,n\in \mathbb{N}$, and $\Phi (l)\in L^{\prime }$ because 
\[
\sum_{k\in \text{supp }l}l(k)\sum_{j_n\in K_n}\delta _{k_{n,0}}(j_n)\delta
_{(k,l,n,j_n)})=\sum_{k\in \text{supp }l}l(k)\delta _{(k,l,n,k_{n,o})} 
\]
and we have assumed that $\delta _{k_{n,0}}\in L_n$, for each $n.$

The next lemma lists some properties of the construction.

\begin{lemma}
\label{Lemma A}Suppose that $K$ and $K_n,n\in \mathbb{N,}$ are compact
Hausdorff spaces, $K_n$ has a distinguished point $k_{n,0}$ and $L,$ $L_n,n%
\mathbb{\in N},$ are sets of purely atomic finitely supported probability
measures on $K$, $K_n$,respectively, with $\delta _{k_{n,0}}\in L_n$, for
each $n,$ as above. Then if $(K^{\prime },L^{\prime })=(K,L)\otimes
\{(K_n,L_n),n\in \mathbb{N\},}$

\begin{enumerate}
\item  $K^{\prime }$ is a compact Hausdorff space and $\phi $ is a
homeomorphism of $K$ into $K^{\prime }$.

\item  \label{Lemma A 2}A net $(k_d,l_d,n_{d,}j_d)_{d\in D}$ in $K^{\prime
}\setminus \Phi (L)$ converges to $(k,l,n,j)$ for some $j\neq k_{n,0}$ if
and only if there exists $d_0\in D$ such that $k_d=k$ , $l_d=l$, and $n=n_d$
for all $d\geq d_0$ and $(j_d)_{d\in D}$ converges to $j$.

\item  A net $(k_d,l_d,n_d,j_d)_{d\in D}$ in $K^{\prime }\setminus \Phi (L)$
converges to $(k,l,n,k_{n,0})=(k,\emptyset )$ if and only if the following
hold.

\begin{enumerate}
\item  $(k_d)_{d\in D}$ converges to $k$ .

\item  If $D_1=\{d:k_d=k\}$ is cofinal in $D$, then for each $l$ and $n$, $%
D_{1,l,n}=\{d\in D_1:l_d=l,n_d=n\}$ is not cofinal in $D$ or $(j_d^{\prime
})_{d\in D}$ converges to $k_{n,0}$, where $j_d^{\prime }=j_d$ if $d\in
D_{1,l,n}$ and $j_d^{\prime }=k_{n,0}$, otherwise.
\end{enumerate}

\item  The map 
\[
l\rightarrow \Phi (l)=\sum_{k\in \text{supp }l}l(k)\delta _{(k,\emptyset )}
\]
for $l\in L$ is a homeomorphism of $L$ into $L^{\prime }$ in the weak$^{*}$
topology.

\item   Each $l^{\prime }\in L^{\prime }$ is atomic and has finite support.

\item  If $L$,$L_n,n\in \mathbb{N}$ are compact in the weak$^{*}$ topology,
then $L^{\prime }$ is compact in the weak$^{*}$ topology.

\item  If $(l_d)$ is a convergent net in $L_n$ with limit $l_0$ and $l\in L$%
, then 
\[
(\sum_{k\in \text{supp }l}l(k)\sum_{j_n\in K_n}l_d(j_n)\delta
_{(k,l,n,j_n)})_d
\]
converges to 
\[
\sum_{k\in \text{supp }l}l(k)\sum_{j_n\in K_n}l_0(j_n)\delta _{(k,l,n,j_n)}
\]
for each $l\in L$.
\end{enumerate}
\end{lemma}

\begin{proof}
We have given a basis for the topology on $K^{\prime }$ in the definition
above. In order to verify the first property we first observe that $%
\{(k,\emptyset ):k\in K\}$ is homeomorphic to $K.$ Notice that the basis for
the topology of $K^{\prime }$ given in the definition above defines the
topology on $\{(k,\emptyset ):k\in K\}$ to be the topology $\{\phi (G):G$ is
open in $K\}$. Thus $\phi $ is a homeomorphism. If $\mathcal{O}\mathrm{\ }$%
is an open cover of $K^{\prime }$ by basic open sets, then there is a finite
subset $\mathcal{O}^{\prime }$ of $\mathcal{O}$ which covers $\phi (K)$. $%
K^{\prime }\setminus \cup \{G_i$, $G_i\in \mathcal{O}^{\prime }\}$, is
contained in a finite union of closed subsets of the form $\{k\}\times
\{l\}\times \{n\}\times F_{k,l,n}$, where $F_{k,l,n}$ is a 
compact subset of $K_n\setminus \{k_{n,0}\}$. The topology on $\{k\}\times
\{l\}\times \{n\}\times F_{k,l,n}$ is the topology induced by identifying
this with $F_{k,l,n}$ in $K_n.$ Therefore a finite number of additional sets
from $\mathcal{O}$ will cover each $\{k\}\times \{l\}\times \{n\}\times
F_{k,l,n}.$ This proves the first assertion.

For the second notice that if $G$ is an open set contained in $K_n\setminus
\{k_{n,0}\}$ and $j\in G$, then $\{k\}\times \{l\}\times \{n\}\times G$ is
an open neighborhood of $(k,l,n,j)$. Thus the net must eventually be in $%
\{k\}\times \{l\}\times \{n\}\times G$. (\ref{Lemma A 2}) follows. Define a
map $\zeta $ from $K^{\prime }$ onto $\{(k,\emptyset ):k\in K\}$ by $\zeta
(k,l,n,j)=k$. Clearly $\zeta $ is continuous and this gives (3)(a), if the
net converges. If $D_1=\{d:k_d=k\}$ is cofinal in $D$, $D_{1,l,n}$ is
cofinal in $D$ and $(j_d^{\prime })_{d\in D}$ had a convergent subnet with
limit $j\neq k_{n,0}$, then there would be an open set $G$ containing $j$
which is contained in $K\setminus \{k_{n,0}\}$. However, $\{k\}\times
\{l\}\times \{n\}\times G$ would be an open set in $K^{\prime }$ containing $%
(k,l,n,j)$ and thus $(k,l,n,j_d^{\prime })_{d\in D}$ would converge to $%
(k,l,n,j)$, which is impossible. Thus (3)(b) holds. Conversely, if we are
given a net satisfying (3)(a) and (b) and $G^{\prime }$ is an open set
containing $(k,l,n,k_{n_0})$, then $G^{\prime }$ contains a
neighborhood of $(k,l,n,k_{n,0})$ of the form 
\begin{eqnarray*}
H &=&\cup _{k^{\prime }\in G}\{(k^{\prime },l^{\prime },n^{\prime
},j_{n^{\prime }}):k^{\prime }\in K,l^{\prime }\in L(k),n^{\prime }\in
\mathbb{N},j_{n^{\prime }}\in K_{n^{\prime }}\} \\
&&\setminus \cup _{(k^{\prime },l^{\prime },n^{\prime })\in F}\{k^{\prime
}\}\times \{l^{\prime }\}\times \{n^{\prime }\}\times F_{k^{\prime
},l^{\prime },n^{\prime }}.
\end{eqnarray*}
Because $(k_d)_{d\in D}$ converges to $k$, there is a $d_0\in D$ such that $%
(k_d,l_d,n_d,j_d)\in H\cup \cup _{(k^{\prime },l^{\prime },n^{\prime })\in
F}\{k^{\prime }\}\times \{l^{\prime }\}\times \{n^{\prime }\}\times
F_{k^{\prime },l^{\prime },n^{\prime }}$ for all $d\geq d_0$. Because $F$ is
finite, we may assume that by choosing another $d_0$ and passing to a subset
of $G$, if necessary, that $F=\{(k,l^{\prime },n^{\prime
}):(l^{\prime },n^{\prime })\in F^{\prime }\}$ for some finite set $%
F^{\prime }$. By (b) we know that for each $(l^{\prime },n^{\prime })\in
F^{\prime }$ there is a $d_{l^{\prime },n^{\prime }}$ such that if $%
(k_d,l_d,n_d,j_d)=(k,l^{\prime },n^{\prime },j_d)$ and $d\geq d_{l^{\prime
},n^{\prime }}$, then $(k_d,l_d,n_d,j_d)\notin \{k\}\times \{l^{\prime
}\}\times \{n^{\prime }\}\times F_{k,l^{\prime },n^{\prime }}$. If $d\geq
d_{l^{\prime },n^{\prime }}$, for all $(l^{\prime },n^{\prime })\in
F^{\prime }$, and $d\geq d_0$, then $(k_d,l_d,n_d,j_d)\in H.$

Because $\phi $ is a homeomorphism (4) is immediate. (5) is obvious from the
definition and the fact that $(k,l,n,j)=(k^{\prime },l^{\prime },n^{\prime
},j^{\prime })$ if and only if $k=k^{\prime }$, $l=l^{\prime }$, $%
n=n^{\prime }$ and $j=j^{\prime }$ or, $j=j^{\prime }=k_{n,0}$ and $%
k=k^{\prime }$. To see that $L^{\prime }$ is compact if $L,L_n,n\in \mathbb{N}$
are. Let $(l_d^{\prime })_{d\in D}$ be a net in $L^{\prime },$ where
\[
l_d^{\prime }=\sum_{k\in \text{supp }l_d}\sum_{j\in
K_{n(d)}}l_d(k)l_{n(d)}^{\prime \prime }(j)\delta _{(k,l_d,n(d).,j)}
\]
for each $d\in D$ . Because $L$ and the $L_n$ are compact, we may assume by
passing to a subnet that the nets $(l_d)_{d\in D}$ and $(l_{n(d)}^{\prime
\prime })_{d\in D}$ converge to $l$ and $l^{\prime \prime }$, respectively.
Here we are thinking of $(l_{n(d)}^{\prime \prime })$ as a net in $\cup
_{n\in \mathbb{N}}\Phi _n(L_n)$, where $\Phi _n$ is the map induced by the
natural embedding $\phi _n$ of $K_n$ into the one point compactification of $%
\cup _{n\in \mathbb{N}}K_n\setminus \{k_{n,0}\}.$ Because $\Phi $ is w$^{*}$%
-continuous, $(\Phi (l_d))_{d\in D}$ converges to $\Phi (l)$. If
$\epsilon>0$ and  $k\in $supp 
$l$, let $G_k$ be an open set containing $k$ such that $l(G_k)<l(k)+\epsilon 
$. We may assume that the sets $G_k$ are disjoint. We must consider two
cases. First suppose that $(l_d)$ has a constant subnet. Then $l_d^{\prime
}=\sum_{k\in \text{supp }l}\sum_{j\in K_{n(d)}}l(k)l_{n(d)}^{\prime \prime}
(j)\delta
_{(k,l_d,n(d),j)}$ for the elements in the subnet and the limit of the
subnet is $\sum_{k\in \text{supp }l}l(k)\sum_{n\in \mathbb{N,}j\in
K_n}l^{\prime \prime }(j)\delta _{(k,l,n.,j)})$. If there is no constant
subnet, then any convergent subnet of points $(k_d,l_d,n_d,j_d)_{d\in D}$
will have a limit of the form $(k,\emptyset )$, where $k$ is the limit of
the first coordinates in the subnet, by (2) and (3). We claim that there is
a convergent subnet with limit $\sum_{k\in \text{supp }l}l(k)\delta
_{(k,\emptyset )}$ Indeed, because $(l_d)$ converges to $l$, there exists a $%
d_0$ such that $l_d(G_k)>l(G_k)-\epsilon $, for all $k\in $supp $l$, $d\geq
d_0.$ This
implies that 
\begin{multline*}
\begin{aligned}
l_d^{\prime }(\cup _{r\in G_k}\{(r,m,n,t) &:r\in K,m\in L(k),n\in
\mathbb{N,} t\in K_n\}\\
&\setminus \cup _{(r,m,n)\in F}\{r\}\times \{m\}\times \{n\}\times
F_{r,m,n})
\end{aligned} \\
>l(G_k)-\epsilon -\sum_{(r,m,n)\in F}l_d^{\prime }(\{r\}\times
\{m\}\times \{n\}\times F_{r,m,n})
\end{multline*}

Notice that $l_d^{\prime }(\{r\}\times \{m\}\times \{n\}\times F_{r,m,n})=0$
if $l_d\neq m$ or $r\notin $supp $l_d$. Because $F$ is finite, and we have
assumed that there is no constant subnet of $(l_d)$ by choosing a another $%
d_1\geq d_0$ we will have $l_d^{\prime }(\{r\}\times \{m\}\times \{n\}\times
F_{r,m,n})=0$ for all $(r,m,n)\in F$ and all $d\geq d_1$. Because $%
l_d^{\prime }$ is a probability measure and $\epsilon >0$ is arbitrary, $%
(l_d^{\prime })$ converges to $\Phi (l)=\sum_{k\in \text{supp }l}l(k)\delta
_{(k,\emptyset )}$. Thus $L^{\prime }$ is compact.

(7) is immediate from the definition. 
\end{proof}

Our next lemma will allow us to compute topological information about the
spaces $K'$ and $L'$ from the component pieces provided the pieces are
properly attached.

\begin{lemma}
\label{AA}Let $K,$ $L,$ $K_{n},$ $L_n,$ $n\in \mathbb{N,}$, $K'$ and $L'$
 be as in the previous
lemma. In addition assume that $K$ is homeomorphic to $[1,\omega ^{\omega
^\alpha m}]$, $K_n$ is homeomorphic to $[1,\omega^{\omega^{\beta(n)}m(n)}]$, $L
$ (with the $\text{w}^*$-topology) 
is homeomorphic to $[1,\omega ^{\omega ^\gamma p}]$,
and $L_n$ is homeomorphic
to $[1,\omega ^{\omega ^{\gamma (n)}p(n)}]$. Moreover, assume that
\begin{align*}
K^{(\omega
^\alpha m)}&=\{k_0\}, \qquad &L^{(\omega ^\gamma p)}&=\{\delta _{k_0}\},\\
K_n^{(\omega^{\beta (n)}m(n))}&=\{k_{n,0}\},\qquad &L_n^{(\omega ^{\gamma
(n)}p(n))}&=\{\delta _{k_{n,0}}\}.
\end{align*}
for all $n.$
 Let 
\[
\omega ^BM=\sup \{\omega ^{\beta
(n)}m(n):n\in \mathbb{N\}} \quad\text{and}\quad
\omega^\Gamma P =\sup \{\omega ^{\gamma
(n)}p(n):n\in \mathbb{N\}}.
\]
Then

\begin{enumerate}
\item  \label{AA1}$K^{\prime (\omega ^BM)}\subset \phi (K)$ and if $\cup \{$%
supp $l:l\in L\}=K$, then $K^{\prime }$ is homeomorphic to $[1,\omega^{\omega
^BM+\omega ^\alpha m}]$ and $K^{\prime (\omega ^BM+\omega ^\alpha m)}=
\{\phi(k_0)\}.$%

\item  \label{AA2}$L^{\prime (\omega ^\Gamma P)}=\Phi (L),$ $L^{\prime }$
 (with the $\text{w}^*$-topology) 
is homeomorphic to $[1,\omega ^{\omega ^\Gamma P+\omega ^\gamma p}]$ and $%
L^{\prime (\omega ^\Gamma P+\omega ^\gamma p)}=\{\delta _{\phi(k_0)}\}$.

\item  \label{AA3}If for each $l\in L$, there is a subset $H_l$ of $K$ such
that $l(H_l)\geq \epsilon $, and $(H_l)_{l\in L}$ are disjoint, and for each 
$n\in \mathbb{N,}l^{\prime \prime }\in L_{n,}$there is a subset $%
H_{n,l^{^{\prime \prime }}}$ of $K_n\setminus \{k_{n,0}\}$ such that $%
l^{\prime \prime }(H_{n,l^{\prime \prime }})\geq \epsilon $, and $%
(H_{n,l^{\prime \prime }})_{l^{\prime \prime }\in L_n},$ are disjoint for
each $n$, then there are disjoint sets $H_{l^{\prime }}^{\prime }$ of $%
K^{\prime }$ for each $l^{\prime }\in L^{\prime }$ such that $l^{\prime
}(H_{l^{\prime }}^{\prime })\geq \epsilon.$ Moreover, if $l\in L$, then we
can define $H_{\Phi (l)}^{\prime }=\phi (H_l)$ and if 
\[
l=\sum_{k\in \text{supp }l}l(k)\sum_{j_n\in K_n}l_n(j_n)\delta
_{(k,l,n,j_n)})
\]
for some $l\in L,n\in \mathbb{N,}l_n\in L_n\setminus \{\delta _{k_{n,0}}\}$,
then we can define
\[
H_{l^{\prime }}^{\prime }=\cup _{k\in \text{supp }l}\{k\}\times \{l\}\times
\{n\}\times H_{l_n}.
\]
\end{enumerate}
\end{lemma}

\begin{proof}
First observe that because $K,L,$ and $K_n$ are countable and the measures
in $L$ are finitely supported, $K^{\prime }$ is countable. If $n\in
\mathbb{N}$, $j\in K_n$ and 
$j\neq k_{n,0}$, then for any $k\in K,l\in L,$ with $l(k)\neq 0$, $%
\{(k,l,n,j^{\prime }):j^{\prime }\neq k_{n,0}\}$ is an open neighborhood of $%
(k,l,n,j)$ in $K^{\prime }$ homeomorphic to $K_n\setminus \{k_{n,0}\}$. Thus 
$(k,l,n,j)$ is in the same derived sets of $K^{\prime }$ as of $K_n$. In
particular, $K_n^{(\omega ^{\beta (n)}m(n))}=\{k_{n,0}\}$ and thus $%
K^{\prime (\omega ^BM)}\subset \phi (K)$. If $\cup \{$supp $l:l\in L\}=K$,
then for each $k\in K$, $\{(k,l,n,j):j\in K_n\}\subset K^{\prime }$ for some 
$l\in L$ and therefore $(k,l,n,k_{n,0})\in K^{\prime (\omega ^\beta M)}$. If 
$k$ is an isolated point in $K$, then $\phi (k)$ is the limit only of
sequences which are eventually in 
\[
\{(k,l,n,j):l\in L,l(k)\neq 0,n\in \mathbb{N,}j\in K_n\}.
\]
Hence $(k,\emptyset )=(k,l,n,k_{n,0})\notin K^{\prime (\omega ^BM+1)}$.
Because $\phi $ is a homeomorphism, it follows that $K^{\prime (\omega
^BM)}=\phi (K^{(0)})$. Similarly, $K^{\prime (\omega ^BM+\rho )}=\phi
(K^{(\rho )})$ for all $\rho $. In particular, $K^{\prime (\omega
^BM+\omega^\alpha m)}=\{\phi(k_0)\}.$

Observe that it follows from Lemma \ref{Lemma A} that $L^{\prime }$ is
countable and compact because $K,K_n$, $L,$and $L_n$ are and thus it is
sufficient to consider the derived sets. If $l\in L$, $n\in \mathbb{N},$ then $%
\{\sum_{k\in \text{supp }l}l(k)\sum_{j\in \text{supp }l_n}l_n(j)\delta
_{(k,l,n,j)}:l_n\in L_n\}$ is homeomorphic (by the obvious map) to $L_n$ and 
\[
\{\sum_{k\in \text{supp }l}l(k)\sum_{j\in \text{supp }l_n}l_n(j)\delta
_{(k,l,n,j)}:l_n\in L_n\}^{(\omega ^{\gamma (n)}p(n))}=\{\Phi (l)\}.
\]
Therefore $\Phi (L)\subset \cap _{n\in \mathbb{N}}L^{\prime (\omega ^{\gamma
(n)}p(n))}=L^{\prime (\omega ^\Gamma P)}$. If $l_n\in L_n\setminus \{\delta
_{k_{n,0}}\}$, then 
\[
\sum_{k\in \text{supp }l}l(k)\sum_{j\in \text{supp }l_n}l_n(j)\delta
_{(k,l,n,j)}-\sum_{k\in \text{supp }l}l(k)l_n(k_{n,0})\delta _{k,\emptyset }
\]
is non-zero and is supported in the open set $\cup _{k\in \text{supp }%
l}\{(k,l,n,j):j\in K_n\setminus \{k_{n,0}\}\}$ and only elements of $%
L^{\prime }$ of the form $\sum_{k\in \text{supp }l}l(k)\sum_{j\in \text{supp 
}l_n^{\prime }}l_n^{\prime }(j)\delta _{(k,l,n,j)}$ with $l_n^{\prime }\in
L_n\setminus \{\delta _{k_{n,0}}\}$ are supported in this set. Therefore $%
\Phi (L)=L^{\prime (\omega ^\Gamma P)}$. Because $\Phi $ is a homeomorphism,
it follows that $\Phi (L)^{(\rho )}=L^{\prime (\omega ^BP+\rho )}$ for all $%
\rho $, proving the second assertion.

For each $l^{\prime }\notin \Phi (L)$, we have defined $H_{l^{\prime
}}^{\prime }$ to be a subset of $K^{\prime }\setminus \phi (K).$ These sets
are clearly disjoint. Also if 
\[
l^{\prime }=\sum_{k\in \text{supp }l}l(k)\sum_{j\in \text{supp }%
l_n}l_n(j)\delta _{(k,l,n,j)}, 
\]
\[
l^{\prime }(H_{l^{\prime }})=\sum_{k\in \text{supp }%
l}l(k)l_n(H_{l_n})=l_n(H_{l_n}). 
\]
Because $\Phi $ is induced by the homeomorphism $\phi $ , $H_{\Phi
(l)}^{\prime }=\phi (H_l),l\in L$, is a family of disjoint subsets of $\phi
(K)$ with $\Phi (l)(H_{\Phi (l)}^{\prime })=l(H_l)$. 
\end{proof}

In order to prove that the evaluation map from $C(K)$ into $C(L)$ is
surjective we will need to show that $L$ is equivalent to the usual unit
vector basis of $l_1.$ The elements of $L$ are not perturbations of
disjointly supported elements and thus the proof uses some special properties of
the construction. We introduce a natural ordering on the elements of $L$
which reflects these properties of the construction.
\begin{definition}
\label{dilate} Suppose $\mathcal{M}$ is a family of measures on a measurable
space $(\Omega ,\mathcal{B)}$ and for each $\mu \in \mathcal{M}$ there is a
set $H_\mu \in \mathcal{B}$ such that $H_\mu \cap H_{\mu ^{\prime
}}=\emptyset $ if $\mu \neq \mu ^{\prime }$, and $\mu (H_\mu )\neq 0.$ Then $%
\mu \succ ^{\prime }\mu ^{\prime }$ if and only if there is a scalar $a\in
(0,|\mu (H_\mu )/(2\mu ^{\prime }(H_{\mu ^{\prime }})|]$ such that $\mu
_{|\cup \{H_{\mu ^{\prime }}:\mu ^{\prime }\neq \mu \}}=a\mu _{|\cup
\{H_{\mu ^{\prime }}:\mu ^{\prime }\neq \mu \}}^{\prime }$ and $|\mu
^{\prime }|(H_\mu )=0$. Define $\mu \succ \nu $ if and only if there is a
finite sequence $(\mu _i)$ in $\mathcal{M}$ such that $\mu =\mu _0\succ
^{\prime }\mu _1\succ ^{\prime }\cdots \succ ^{\prime }\mu _k=\nu $ $.$
\end{definition}

Notice that $\mu \succ \mu $ is impossible and the relation $\succ $ is
transitive by definition. Thus we can define a partial order on $\mathcal{M}$
by $\mu \succeq \mu ^{\prime }$ if and only if $\mu =\mu ^{\prime }$ or $\mu
\succ \mu ^{\prime }.$ Although the relation is really on the pairs
$(\mu,H_\mu)$, we will write it as though it were on the measures. This
will not present any difficulty because the sets $H_\mu$ will be fixed
during the construction.

The relation above occurs naturally in the construction of the pairs $(K,L).$
For $(K^{\prime },L^{\prime })=(K,L)\otimes \{(K_n,L_n),n\in \mathbb{N\}}$ as
in Lemma \ref{Lemma A} each $l^{\prime }\in L^{\prime }$ which is of the
form $\sum_{k\in \text{supp }l}l(k)\sum_{j_n\in K_n}l_n(j_n)\delta
_{(k,l,n,j_n)})$ for some $l\in L,n\in \mathbb{N,}l_n\in L_n$, satisfies $%
l_{|K^{\prime }\setminus (\text{supp }l)\times \{l\}\times \{n\}\times
K_n}^{\prime }=l_n(k_{n,0})l$. If we have the sets $(H_{l^{\prime }}^{\prime
})_{l^{\prime }\in L^{\prime }}$, defined as in Lemma \ref{AA}, (\ref{AA3}),
and for $l^{\prime \prime }\in L_n$, supp $l^{\prime \prime }\subseteq
H_{n,l^{\prime \prime }}\cup \{k_{n,0}\}$ then $l^{\prime }\succ ^{\prime }l.
$

The next lemma is similar to Proposition IV.13 of \cite{Gasp}. It will be
used to show that the sets of measures $L$ that we construct actually are
equivalent to the basis of $l_1.$

\begin{lemma}
\label{B}Suppose that $M$ is a set of mutually singular probability measures
on a measurable space $(K,\mathcal{B)},$ $\epsilon >0,$ and that $(\mu _n)$
is a sequence of (finite) convex combinations of the measures in $M$ and ($%
A_n)$ is a sequence of disjoint measurable sets. Let $\succ ^{\prime }$ and $%
\succ $ be defined as above for $\mathcal{M=}\{\mu _n\}$ and $H_{\mu _n}=A_n.
$ Suppose that $(\mu_n,A_n)_{n=1}^\infty$ satisfy the following.

\begin{enumerate}
\item  \label{B1} For each $n\in \mathbb{N},$ $\mu _n(A_n)\geq \epsilon .$

\item  \label{B2} For each $n\in \mathbb{N},$ either there is a unique $%
n^{\prime }\in \mathbb{N}$ such that $\mu _n\succ ^{\prime }\mu _{n^{\prime }}$
or for all $n^{\prime }\neq n,\mu _n(A_{n^{\prime }})=0.$

\item  \label{B3} For all $n\neq m$ if it is not the case that $\mu _n\succ \mu
_m$, then $\mu _n(A_m)=0.$
\end{enumerate}

Then $||\sum c_n\mu _n||\geq (2\epsilon /3)\sum |c_n|$ for any sequence of
scalars $(c_n).$
\end{lemma}

\begin{proof}
Because we only use information about the measures on the sets $A_n$,
without loss of generality we may assume that $\mu _n\in 
\text{co}\{m\in \mathcal{%
M}:m(A_k)>0,$ for some $k\}\cup \{0\}.$ The relation $\succ $ defines a
partial order on $\{\mu _n\}$. Observe that if $\mu \succ ^{\prime }\nu $
and $\mu =\sum_{j\in F}b_j^\mu m_j$ and $\nu =\sum_{j\in G}b_j^\nu m_j$,
where $m_j\in \mathcal{M}$ and $b_j^\mu ,b_j^\nu $ are non-zero for all $j$,
then $F\supset G$ Therefore, since each $\mu _n$ is a finite convex
combination, for any $n(0)\in \mathbb{N}$ there is a unique finite maximal
sequence $(\mu _{n(i)})_{i=0}^k$ such that $\mu _{n(0)}\succ ^{\prime }\mu
_{n(1)}\succ ^{\prime }\cdots \succ ^{\prime }\mu _{n(k)}$. Let $(c_n)$ be a
finite sequence of scalars and let 
\[
F=\{n:\exists n^{\prime \text{ }}\text{ such that }c_{n^{\prime }}\neq 0%
\text{ and }\mu _{n^{\prime }}\succ \mu _n\}.
\]
Clearly $F$ is a finite set. Partition $F$ into sets $(F_j)_{j=0}^J$ such
that for each $j<J$ and $n\in F_j$ there is an $n^{\prime }\in F_{j+1}$ such
that $\mu _n\succ ^{\prime }\mu _{n^{\prime }}$ and for all $n^{\prime }\neq
n,n^{\prime }\in F_j,$ $\mu _{n^{\prime }}$ and $\mu _n$ are incomparable.
If $\mu _n\succ ^{\prime }\mu _{n^{\prime }}$, let $a_{n,n^{\prime }}$
denote the scalar such that $\mu _{n|\cup \{A_k:k\neq n\}}=a_{n,n^{\prime
}}\mu _{n^{\prime }|\cup \{A_k:k\neq n\}}$. For notational convenience, let
$a_{n,n^{\prime }}=0$ if it is not the case that  $\mu _n\succ ^{\prime
}\mu _{n^{\prime }.}$
A simple induction argument
using \ref{B2} and \ref{B3} shows that
\begin{align*}
\|\sum c_n\mu _n\| &= \|\sum_{j=0}^J\sum_{n(j)\in
F_j}c_{n(j)}\mu_{n(j)}\| \\
&= \|\sum_{j=0}^J\sum_{n(j)\in
F_j}c_{n(j)}{\mu_{n(j)}}_{|\cup \{A_n:\mu_{n(j)}\succeq \mu_n\}}\| \\
&=\|\sum_{j=0}^J\sum_{n(j)\in
F_j}\sum_{n:\mu_n\succeq \mu_{n(j)}} c_{n}{\mu_{n}}_{|A_{n(j)}}\|.
\end{align*}
Another induction argument and the definition of the scalars
$a_{n,n^{\prime }}=0$ gives the following inequality.
\begin{eqnarray*}
\|\sum c_n\mu _n\| \geq \sum_{j=0}^J\sum_{n(j)\in
F_j}|c_{n(j)}&+&\sum_{n(j-1)\in F_{j-1}}a_{n(j-1),n(j)}(c_{n(j-1)} \\
&+&\sum_{n(j-2)\in F_{j-2}}a_{n(j-2),n(j-1)}(c_{n(j-2)}+\cdots  \\
&+&\sum_{n(0)\in F_0}a_{n(0),n(1)}c_{n(0)}))|\mu _{n(j)}(A_{n(j)}).
\end{eqnarray*}

Now we split off $1/3$ of each term and shift the index on these pieces 
to combine with the subsequent related term.

\begin{multline*}
\begin{aligned}
\sum_{j=0}^J\sum_{n(j)\in F_j} \biggl | c_{n(j)} &+ \sum_{n(j-1)\in
F_{j-1}}a_{n(j-1),n(j)}(c_{n(j-1)} 
\\ &+\sum_{n(j-2)\in
F_{j-2}}a_{n(j-2),n(j-1)}(c_{n(j-2)}+\cdots
 \\  &+\sum_{n(0)\in
F_0}a_{n(0),n(1)}c_{n(0)})) \biggr| \mu _{n(j)}(A_{n(j)}) 
\end{aligned}
 \\
\begin{aligned} =\sum_{j=0}^J\sum_{n(j)\in F_j}\biggl (\frac 23 \biggl | c_{n(j)}
&+\sum_{n(j-1)\in F_{j-1}}a_{n(j-1),n(j)}(c_{n(j-1)} \\
&+\sum_{n(j-2)\in
F_{j-2}}a_{n(j-2),n(j-1)}(c_{n(j-2)}+\cdots 
\\ &+\sum_{n(0)\in
F_0}a_{n(0),n(1)}c_{n(0)}))\biggr|\mu _{n(j)}(A_{n(j)}) \\ &+\frac 13\sum_{\mu
_{n(j-1)}\succ ^{\prime } \mu_{n(j)}}\biggl|c_{n(j-1)}\\ &+\sum_{n(j-2)\in
F_{j-2}}a_{n(j-2),n(j-1)}(c_{n(j-2)}
\\ &+\sum_{n(j-3)\in
F_{j-3}}a_{n(j-3),n(j-2)}(c_{n(j-3)}+\cdots \\ &+\sum_{n(0)\in
F_0}a_{n(0),n(1)}c_{n(0)})) \biggr |
\mu _{n(j-1)}(A_{n(j-1)})\biggr ).
\end{aligned}
\end{multline*}

The condition $\mu _{n(j-1)}\succ ^{\prime }\mu _{n(j)}$ is equivalent to $%
a_{n(j-1),n(j)}\neq 0$ and by the definition of $\succ ^{\prime }$, $%
2a_{n(j-1),n(j)}\mu _{n(j)}(A_{n(j)})\leq \mu _{n(j-1)}(A_{n(j-1)})$.
Therefore by the triangle inequality, 
\[
\|\sum c_n\mu _n\|\geq \sum_{j=0}^J\sum_{n(j)\in F_j}\frac 23|c_{n(j)}|\mu
_{n(j)}(A_{n(j)})\geq \sum_{j=0}^J\sum_{n(j)\in F_j}\frac
23|c_{n(j)}|\epsilon . 
\]

\end{proof}

\section{Construction of the Operators}

The aim of this section is to produce pairs $%
(K_\alpha ,L_\alpha )_{\alpha <\omega _1}$ by transfinite induction
such that $K_\alpha $ is a
countable compact Hausdorff space and $L_\alpha $ is a w$^{*}$-closed subset
of the probability measures in $C(K_\alpha )^{*}$ which is equivalent to the
basis of $l_1.$ 

Fix an ordinal $\zeta <\omega _1.$ Let $\zeta _n\uparrow \omega ^\zeta $ and
for each $n\in \mathbb{N}$ let $S_n=[1,\omega ^{\zeta _n}]$ with the order
topology and let $T_n=\{\frac 12(\delta _\beta +\delta _{\omega ^{\zeta
_n}}):\beta \leq \omega ^{\zeta _n}\}$. Let the distinguished point of $S_n$
be $\omega ^{\zeta _n}$. Let $S_0=[1,1]$ and $T_0=\{\delta _1\}$. Define 
\[(K_1,L_1)=(S_0,T_0)\bigotimes \{(S_n,T_n):n\in \mathbb{N}\}.
\]
It is easy to see
that up to a homeomorphism of $[1,\omega ^{\omega ^\zeta }]$ we could have
defined $K_1=[1,\omega ^{\omega ^\zeta }]$ and $L_1=\{\frac 12(\delta _\beta
+\delta _{\omega ^{\omega ^\zeta }}):\beta \leq \omega ^{\omega ^\zeta }\}$.
We take the distinguished point $k_1$ of $K_1$ to be $\phi (1)$ where $1\in
S_0$. Now suppose that we have defined $K_\gamma $ and $L_\gamma $ for all $%
\gamma <\alpha .$ Let $k_\gamma $ denote the distinguished point of $%
K_\gamma $.There are two cases. First assume that $\alpha =\alpha ^{\prime
}+1$ for some $\alpha ^{\prime }$. Define  
\[
(K_\alpha ,L_\alpha
)=(K_1,L_1)\bigotimes (K_{\alpha ^{\prime }},L_{\alpha ^{\prime }}). 
\]
Let the distinguished point be  $k_\alpha =\phi (k_1)$
(More formally we should have a sequence of spaces $%
\{(K_{\alpha _n},L_{\alpha _n}):n\in \mathbb{N}\}$ on the right of $%
\bigotimes $, but we can take $(K_{\alpha _1},L_{\alpha _1})=$ $(K_{\alpha
^{\prime }},L_{\alpha ^{\prime }})$ and $(K_{\alpha _n},L_{\alpha
_n})=(S_0,T_0)$ for $n>1$. These spaces $(S_0,T_0)$ have no effect on $%
(K_\alpha ,L_\alpha )$.) If $\alpha $ is a limit ordinal, let $(\alpha _n)$
be an increasing sequence of ordinals with limit $\alpha $. Let 
\[
(K_\alpha
,L_\alpha )=(S_0,T_0)\bigotimes \{(K_{\alpha _n},L_{\alpha _n}):n\in \mathbb{N}
\}.
\]
 Let $\phi (1),$ for $1\in S_0,$ be the distinguished point of $K_\alpha 
$. The definition for $\alpha$ a limit ordinal depends on the sequence
$(\alpha_n)$. However, the properties of the space are not dependent on the
sequence and we will assume that whenever we use a sequence approaching
$\alpha$ that it is the same one. 
This completes the definition of the pairs $(K_\alpha ,L_\alpha )$. Notice
that we actually have such a transfinite family of spaces for each $\zeta
<\omega _1$. The choice of $\zeta $ will be made in the proof of Theorem \ref
{main}

Now we must consider the properties of these pairs. First we compute the
topological information using Lemma \ref{AA}. As noted above $K_1$ and $L_1$
are homeomorphic to $[1,\omega ^{\omega ^\zeta }]$. Notice that we have the
following relations. If $K_{\alpha ^{\prime }}$ and $L_{\alpha ^{\prime }}$
are homeomorphic to $[1,\omega ^{\omega ^\zeta \beta }]$, then $K_{\alpha
^{\prime }+1}$ and $L_{\alpha ^{\prime }+1}$ are homeomorphic to $[1,\omega
^{\omega ^\zeta (\beta +1)}]$, by Lemma \ref{AA}, (\ref{AA1}) and (\ref{AA2}%
). If $(\alpha _n)$ is an increasing sequence of ordinals with limit $\alpha 
$ and $K_{\alpha _n}$ and $L_{\alpha _n}$ are homeomorphic to $[1.\omega
^{\omega ^\zeta \beta _n}]$, then $K_\alpha $ and $L_\alpha $ are
homeomorphic to $[1,\omega ^{\omega ^\zeta \beta }]$ where $\beta =\sup
\beta _n$. Therefore a straight-forward transfinite induction argument shows
that $K_\alpha $ and $L_\alpha $ are homeomorphic to $[1,\omega ^{\omega
^\zeta \alpha }]$ for all $\alpha <\omega _1$.

Next we will show that $L_\alpha $ is equivalent to the standard basis of $%
l_1.$ We will use Lemma \ref{AA}, (\ref{AA3}) and Lemma \ref{B}. For $\delta
_1\in T_0$ we take $H_{\delta _1}=S_0.$ For each $n\in \mathbb{N},$ $\frac
12(\delta _\beta +\delta _{\omega ^{\zeta _n}})\in T_n$ and $\beta \leq
\omega ^{\zeta _n}$ , we let $H_{\frac 12(\delta _\beta +\delta _{\omega
^{\zeta _n}})}=\{\beta \}.$ If $l\in L_1$ , then $H_l^1=\phi (S_0)$ if $%
l=\Phi (\delta _1)=\delta _{(1,\emptyset )},$ and $H_l^1=\{_{(1,\delta
_1,n,\beta )}\}$ if $l=\frac 12(\delta _{(1,\delta _1,n,\beta )}+\delta
_{(1,\emptyset )}),$ as in Lemma \ref{AA}. Notice that for all $l\neq \delta
_{(1,\emptyset )},$ $l\in L_1$, $l\succ ^{\prime }\delta _{(1,\emptyset )}.$
Thus by Lemma \ref{B} with $(\mu _n)=(\delta _{(1,\delta _1,n,\beta
)})_{\beta \leq \omega ^{\zeta _n},n\in \mathbb{N}},$ $L_1$ is $3$%
-equivalent to the basis of $l_1.$ Notice that $\Phi (\delta _1)$ is the
only element of $L_1$ which does not have a successor (under $\succ ^{\prime
})$ and that if $l\succ ^{\prime }l^{\prime },$ then $l_{|\cup \{H_m:m\neq
l\}}=\frac 12l^{\prime }.$ 

Assume inductively that for each $\beta <\alpha ,$ we have defined sets $%
H_l^\beta \subset K_\beta $ for all $l\in L_\beta $, satisfying the
hypothesis of Lemma \ref{B} with $\epsilon =\frac 12,$ $(\mu _n)$ as the
point mass measures on $K_\beta $ and $l\succeq \delta _{k_\beta },$ for all 
$l\in L_\beta $. Further assume that if $l,l'\in L_\beta$ and
$l\succ ^{\prime }l^{\prime },$ then 
$l_{|\cup \{H_m:m\neq l\}}=\frac 12l^{\prime }.$ Define the sets $%
H_{l^{\prime }}^\alpha ,$ $l^{\prime }\in L_a,$ as in Lemma \ref{AA}, (\ref
{AA3}). Thus Lemma \ref{B} (\ref{B1}) is satisfied. We need to verify the
other hypotheses of Lemma \ref{B}. In order to handle the successor ordinal
case and the limit ordinal case, at the same time, let $O_1=K_{\alpha
^{\prime }}$ and $P_1=L_{\alpha ^{\prime }},$ and $O_n=S_0$ and $P_n=T_0$
for all $n>1,n\in \mathbb{N}$, in the case $\alpha =\alpha ^{\prime }+1,$
and let $O_n=K_{\alpha _n}$ and $P_n=L_{\alpha _n}$., if $\alpha =\lim $ $%
\alpha _n.$ Also let $O_0=K_1,$ $P_0=L_1,$ or $O_0=S_0,$ $P_0=T_0,$
respectively. Let $o_n$ be the distinguished point of $O_n$ for $%
n=0,1,\ldots $ .With this notation $(K_\alpha ,L_\alpha
)=(O_0,P_0)\bigotimes \{(O_n,P_n):n\in \mathbb{N}\}.$ For each $l\in P_n$ let
$H_l^n$ be the associated subset of $O_n.$

Because $\phi $ is a homeomorphism, it follows that $\Phi (O_0)$ and $%
H_l^0,l\in \Phi (O_0),$ satisfy the hypothesis of Lemma \ref{B}. Moreover $%
\phi (l_1)\succ ^{\prime }\phi (l_2)$ if and only if $l_1\succ ^{\prime }l_2.
$ Suppose $l_0,l_0^{^{\prime \prime }}\in O_0$ , $l_n\in P_n,l_{n^{\prime
\prime }}^{\prime \prime }\in P_{n^{\prime \prime }}$ , for some $%
n,n^{\prime \prime }.$ and 
\[
l^{\prime }=\sum_{k\in \text{supp }l_0}l_0(k)\sum_{j_n\in K_n}l_n(j_n)\delta
_{(k,l_0,n,j_n)})
\]
and 
\[
l^{\prime \prime }=\sum_{k\in \text{supp }l_0^{\prime \prime }}l_0^{\prime
\prime }(k)\sum_{j_{n}\in K_{n^{\prime \prime
}}}l_{n^{\prime \prime }}^{\prime \prime }(j_{n})\delta
_{(k,l_0^{\prime \prime },n^{\prime \prime },j_{n})}.
\]
Then $l^{\prime }\succ l_0$ and $l^{\prime \prime }\succ l_0^{\prime \prime
}.$ If $l_0\neq l_0^{\prime \prime }$ or $n\neq n^{\prime \prime }$ , then $%
H^\alpha_{l^{\prime }}\cap
H^{\alpha}_{l^{\prime \prime }}=\emptyset $ because $H^\alpha_{l^{\prime
}}=\cup _{k\in \text{supp }l_0}\{k\}\times \{l_0\}\times \{n\}\times H^n_{l_n}$
and $H^{\alpha}_{l^{\prime \prime }}=\cup _{k\in \text{supp }l_0^{\prime \prime
}}\{k\}\times \{l_0^{\prime \prime }\}\times \{n^{\prime \prime }\}\times
H^{n''}_{l_{n''}^{\prime \prime }}.$
Also $l^{\prime }(H^\alpha_{l^{\prime \prime }})=0$ and $%
l^{\prime \prime }(H^\alpha_{l^{\prime }})=0$.
If $l_0=l_0^{\prime \prime }$ and $%
n=n^{\prime \prime },$ but $l_n\neq l_{n^{\prime \prime }}$ then $l^{\prime
}\succ l^{\prime \prime }$ if and only if $l_n\succ l_{n^{\prime \prime }}.$
Moreover, because $P_n$ and $(H^n_l)_{l\in P_n}$ satisfy the hypotheses of
Lemma \ref{B}, 
\[
l_0\bigotimes P_n=\{\sum_{k\in \text{supp }l_0}l_0(k)\sum_{j_n\in
K_n}l_n(j_n)\delta _{(k,l_0,n,j_n)}):l_n\in P_n\}
\]
and $(H_{l^{\prime }})_{l^{\prime }\in l_0\bigotimes P_n}$ satisfy the same
conditions and 
\[
\sum_{k\in \text{supp }l_0}l_0(k)\sum_{j_n\in K_n}l_n(j_n)\delta
_{(k,l_0,n,j_n)}\succ l_0=\sum_{k\in \text{supp }l_0}l_0(k)\sum_{j_n\in
K_n}\delta _{(o_n\}}(j_n)\delta _{(k,l_0,n,j_n)}
\]
for all $l_n.$ Observe that Lemma \ref{B} (\ref{B2}) is therefore satisfied.
Indeed, if $l_n\neq \delta _{o_n}$ , then there is some $l_n^{\prime }\in P_n
$ such that $l_n\succ ^{\prime }l_n^{\prime }$ and thus 
\[
\sum_{k\in \text{supp }l_0}l_0(k)\sum_{j_n\in K_n}l_n(j_n)\delta
_{(k,l_0,n,j_n)}\succ ^{\prime }\sum_{k\in \text{supp }l_0}l_0(k)\sum_{j_n%
\in K_n}l_n^{\prime }(j_n)\delta _{(k,l_0,n,j_n)}.
\]
If $l_n=\delta _{o_n},$ then $l_0=\sum_{k\in \text{supp }l_0}l_0(k)\sum_{j_n%
\in K_n}\delta _{(o_n\}}(j_n)\delta _{(k,l_0,n,j_n)})$ and there is some $%
l_0^{\prime }\in O_0$ such that $l_0\succ ^{\prime }l_0^{\prime }$ or $%
l_0=\delta _{o_0}.$ Further by transitivity of $\succ $ it follows that for
all $l\in L_\alpha ,$ $l\succeq \Phi (\delta _{o_0}).$ For (\ref{B3}) we
need only consider the case of an element of the form $\Phi (l)$ for some $%
l\in P_0$ and an element of the form $l^{\prime }=\sum_{k\in \text{supp }%
l_0}l_0(k)\sum_{j_n\in K_n}l_n(j_n)\delta _{(k,l_0,n,j_n)}$ where $l_0\in P_0
$ and $l_n\in P_n$ for some $n>1.$ In this case there are three
possibilities (a) $l\succ l_0,$ (b) $l_0\succeq l,$ or (c) neither (a) nor
(b). Case (b) gives $\sum_{k\in \text{supp }l_0}l_0(k)\sum_{j_n\in
K_n}l_n(j_n)\delta _{(k,l_0,n,j_n)}\succeq \Phi (l)$ and so there is nothing
to do. In case (c) supp $l^{\prime }\subseteq \cup _{k\in \text{supp }%
l_0}\{k\}\times \{l_0\}\times \{n\}\times O_n\cup $ supp $\Phi (l_0)$ and $%
H_{\Phi (l)}\subseteq \phi (O_0)\setminus $ supp $\Phi (l_0).$ Therefore $%
\Phi (l)(H_{l^{\prime }})=0$ and $l^{\prime }(H_{\Phi (l)})=0.$ In case (a) $%
\alpha =\alpha ^{\prime }+1$ and thus $l_0=\delta _{(1,\emptyset )}$ and $%
l=\frac 12(\delta _o+\delta _{(1,\emptyset )})$ for some $o\in K_1\setminus
\{(1,\emptyset )\}.$ Clearly $\Phi (l)(H_{l^{\prime }})=0$ and $l^{\prime
}(H_{\Phi (l)})=0$ in this case also.

We have thus proved the following.

\begin{proposition}
\label{l1} For each $\zeta <\omega _1$ there is a family of pairs $(K_\alpha
,L_\alpha )_{\alpha <\omega _1}$ , where for each $\alpha $ $K_\alpha $ is
homeomorphic to $[1,\omega ^{\omega ^\zeta \alpha }]$ and $L_\alpha $ is a w$%
^{*}$ -closed subset of the probability measures in $C(K_\alpha )^{*}$ which
is homeomorphic to $[1,\omega ^{\omega ^\zeta \alpha }]$ in the w$^{*}$
-topology. Moreover, $L_\alpha $ is $3$-equivalent to the usual basis of $%
l_1.$ Consequently, the evaluation map $T:C(K_\alpha )\rightarrow C(L_\alpha
)$ defined by $T(f)(l)=l(f)$ , for all $l\in L_\alpha ,$ is a surjection.
\end{proposition}

\begin{remark}
We actually have that $[L_\alpha ]$ is isometric to $l_1.$ To see this
observe that for $\alpha =1$ the elements $(l_{|H_l})_{l\in L_1}$ are
disjointly supported elements of $[L_1]$ with span containing $L_1.$ Thus
the normalized sequence is a basis for $[L_1]$ which is $1$-equivalent to
the basis of $l_1.$ An induction argument shows that for all $\alpha <\omega
_1,$ $(l_{|H_l})_{l\in L_\alpha }$ are disjointly supported elements of $%
[L_\alpha ]$ with span containing $L_\alpha $ and thus $[L_\alpha ]$ is
isometric to $l_1.$ Notice that this also means that the argument about the
equivalence of $L_\alpha $ to the usual $l_1$ basis could have been made
using $(l_{|H_l})_{l\in L_\alpha }$ in the role of $(\mu _n)$ in the
application of Lemma \ref{B} rather than using the point mass measures.
\end{remark}

\section{The Wolfe index of operators}

In this section we will show that the evaluation operators defined in
Proposition \ref{l1} are actually small in the sense that for most $\alpha $
the ordinals $\beta $ for which there is a subspace $X$ of $C(K_\alpha )$
which isomorphic to $C(\omega ^\beta )$ and for which $T_{|X}$ is an
isomorphism are much smaller than $\omega ^\zeta \alpha .$ The device for
computing the possible ordinals $\beta $ is an ordinal index which was
defined in \cite
{wolfe} and characterized in \cite{A1}.

\begin{definition}
\label{wi} Let $K$ be a compact Hausdorff space, $\epsilon >0,$ and let $B$
be a subset of $C(K)^{*}.$ Let 
\[
P_0(\epsilon ,B)=\{(\mu ,G):\mu \in B,G\text{ is open in }K,|\mu |(G)\geq
\epsilon \}.
\]
If $P_\alpha (\epsilon ,B)$ has been defined, let 
\begin{equation*}
\begin{split}
P_{\alpha +1}(\epsilon ,B) &=\{(\mu ,G)\in P_0(\epsilon ,B)\text{ and there
is a sequence } \\
&(\mu _n,G_n)_{n=1}^\infty  \subset P_\alpha (\epsilon ,B)\text{ such that }%
\mu _n\stackrel{\text{w}^{*}}{\rightarrow }\mu , \\
&G_n\cap G_{n^{\prime }} =\emptyset ,\text{ for }n\neq n^{\prime },\text{
and }\overline{\cup G_n}\subset G\}.
\end{split}
\end{equation*}
For a limit ordinal $\beta $ let 
\begin{equation*}
\begin{split}
P_\beta (\epsilon ,B) &=\{(\mu ,G)\in P_0(\epsilon ,B)\text{ and there is a
sequence of ordinals }\\ 
&\alpha _n\uparrow \beta \text{ and} 
(\mu _n,G_n) \subset P_{\alpha _n}(\epsilon ,B)\text{ such that }\\
&\mu _n%
\stackrel{\text{w}^{*}}{\rightarrow }\mu , 
G_n\cap G_{n^{\prime }} =\emptyset ,\text{ for }n\neq n^{\prime },\text{
and }\overline{\cup G_n}\subset G\}.
\end{split}
\end{equation*}

\end{definition}

The result that we will use here is the following.

\begin{theorem}
\label{AW} Let $T$ be a bounded operator from $C(K)$ into a separable Banach
space $X.$ Then there is a subspace $Y$ of $C(K)$ such that $Y$ is
isomorphic to $C(\omega ^{\omega ^\alpha })$ and $T_{|Y}$ is an isomorphism
if and only if there is an $\epsilon >0$ such that $P_\gamma (\epsilon,
T^{*}(B_{X^{*}})\neq \emptyset $ for all $\gamma
<\omega ^\alpha .$
\end{theorem}

This result is an amalgamation of Theorems 0.2 and 0.3 from \cite{A1}. It
follows that we need only bound the Wolfe index. As in the previous
section we keep $\zeta $ fixed and consider the evaluation operators $%
T_\alpha :C(K_\alpha )\rightarrow C(L_\alpha ).$ In this case the expression 
$T^{*}(B_{X^{*}})$ which occurs in Theorem \ref{AW}
is $T_\alpha ^{*}$ of the unit ball of $C(L_\alpha )^{*},$ which is w$^{*}$
-closed. Also if $\mu \in B_{C(L_\alpha )^{*}},$ then $\mu =\sum_{l\in
L_\alpha }c_l\delta _l,$ where $\sum_{l\in L_\alpha }|c_l|\leq 1.$ Hence $%
T^{*}(\mu )=\sum_{l\in L_\alpha }c_ll.$ These observations will allow us to
employ the following lemma from \cite{A1} (Lemma 3.2) to reduce to
considering only the sets $L_\alpha $.

\begin{lemma}
\label{coL} Let $L$ be a w$^{*}$ -closed countable subset of $\{\mu :\mu \in
B_{C(K)^{*}},\mu >0\}$ for some countable compact metric space $K.$ Suppose
that the evaluation map $T:C(K)\rightarrow C(L)$ defined by $(Tf)(l)=l(f)$ ,
for all $l\in L,$ is surjective. Then, for $\alpha <\omega _1,$ there is an $%
\epsilon >0,$ such that $P_\gamma (\epsilon ,\overline{\text{co}}(\pm
L))\neq \emptyset ,$ for all $\gamma <\omega ^\alpha ,$ if and only if there
is an $\epsilon ^{\prime }>0,$ such that $P_\gamma (\epsilon ^{\prime
},L)\neq \emptyset ,$ for all $\gamma <\omega ^\alpha .$
\end{lemma}

Before we apply this to the examples let us make a few observations about
the sets $P_\gamma (\epsilon ,L_\alpha ).$ Because the sequence $(\mu _n)$
occurring in the definition is a sequence of distinct elements, for any
ordinals $\gamma $ and $\eta ,$ 
\[
\{l:(l,G)\in P_{\gamma +\eta }(\epsilon ,L_\alpha )\}\subseteq \{l:(l,G)\in
P_\gamma (\epsilon ,L_\alpha )\}^{(\eta )}.
\]
Also, because the sets $P_\gamma (\epsilon ,L_\alpha )$ decrease to $%
\emptyset ,$ 
\[
\{l:(l,G)\in P_\gamma (\epsilon ,L_\alpha )\}\setminus
\{l:(l,G)\in P_{\gamma +1}(\epsilon ,L_\alpha )\}
\]
 is dense in $\{l:(l,G)\in
P_\gamma (\epsilon ,L_\alpha )\}$ for all $\gamma .$ Moreover, if $(l,G)\in
P_{\gamma +1}(\epsilon ,L_\alpha ),$ there exists $((l_n,G_n))_{n=1}^\infty
\subseteq P_\gamma (\epsilon ,L_\alpha )\setminus P_{\gamma +1}(\epsilon
,L_\alpha )$ such that $l_n\stackrel{\text{w}^{*}}{\rightarrow }%
l,$ $l_n(G_n)\geq \epsilon,$ $l(G)\geq \epsilon,$ $G_n\cap G_{n^{\prime
}}=\emptyset ,$ for $n\neq n^{\prime },$ and $\overline{\cup G_n}\subset G.$

\begin{definition}
For each $\epsilon >0$ and $\alpha <\omega _1$ let $\rho (\epsilon ,\alpha
)=\sup \{\gamma :P_\gamma (\epsilon ,L_\alpha )\neq \emptyset \}.$
\end{definition}

It is easy to see that for $\alpha =1$ and $\frac 12<\epsilon ,$ $\rho
(\epsilon ,1)=0$ and for $0<\epsilon \leq \frac 12,$ $\rho (\epsilon
,1)=\omega ^\zeta .$ Obviously $\rho (\epsilon ,\alpha )=0$ for all $\alpha $
if $\epsilon >1.$ The lemma below will permit us to estimate $\rho (\epsilon
,\alpha )$ for all $\alpha <\omega _1$ and $\epsilon \leq 1.$

\begin{lemma}
\label{ind} For each $\beta <\omega _1,$ $\rho (\epsilon ,\beta +1)\leq \max
(\rho (2\epsilon ,\beta )+\rho (\epsilon ,1),\rho (\epsilon ,\beta )+1).$ If 
$\beta _n\uparrow \beta ,$ $\rho (\epsilon ,\beta )=\lim \rho (\epsilon
,\beta _n).$
\end{lemma}
\begin{proof}
Suppose that $\alpha =\beta +1$ for some $\beta <\omega _1.$ Before we begin
estimating $\rho$,
 let us look at relationship between $P_\gamma (\epsilon ,L_\alpha
\setminus \Phi (L_1))$ and pairs $(l,G)$ with $l\in \Phi (L_1).$

Let $(l_n)_{n\in \mathbb{N}}\subset L_\alpha
\setminus \Phi (L_1)^{(1)},l\in \Phi
(L_1)$ and let $(G_n)_{n\in \mathbb{N}}$ and $G$ be
open subsets of $K_\alpha $ such
that $l_n\stackrel{\text{w}^{*}}{\rightarrow }l,l_n(G_n)\geq \epsilon
,l(G)\geq \epsilon ,G_n\cap G_{n^{\prime }}=\emptyset ,$ for $n\neq
n^{\prime },$ and $\overline{\cup G_n}\subset G.$ Then for some $k\in K$, $%
l=\frac 12(\delta _k+\delta _{k_1}),$ and by passing to a subsequence we may
assume that 
\[
l_n=\frac 12(\sum_{j\in K_\beta }l_n^{\prime }(j)\delta
_{(k(n),l_n^{^{\prime \prime }},m,j)}+\sum_{j\in K_\beta }l_n^{\prime
}(j)\delta _{(k_1,l_n^{^{\prime \prime }},m,j)}),
\]
 for some $l_n^{\prime
}\in L_\beta ,l_n^{^{\prime \prime }}=\frac
12(\delta _{k(n)}+\delta _{k_1})$ with $l_n^{^{\prime \prime }}\stackrel{%
\text{w}^{*}}{\rightarrow }l.$ (Because of the definition of the pair $%
(K_\alpha ,L_\alpha )$ for $\alpha$ a successor ordinal,
only the value $1$ of the third index $m$ is of any
interest.) For each $n$ let
\begin{align*}
l_n^1&=\sum_{j\in K_\beta }l_n^{\prime
}(j)\delta _{(k(n),l_n^{\prime \prime },m,j)} \\
\intertext{ and } 
l_n^2&=\sum_{j\in K_\beta }l_n^{\prime }(j)\delta _{(k_1,l_n^{\prime \prime
},m,j)}.
\end{align*}
Then $l_n^1(G_n)\geq \epsilon$ or $l_n^2(G_n)\geq \epsilon$ for each $n$
and thus for one of $(l^1_n)$ and $(l_n^2)$ there are infinitely many
such $n.$

Suppose that for some $\gamma<\omega_1,$
$(l_n,G_n)\in P_\gamma(\epsilon,L_\alpha)$ for each $n$.
Because $L_\alpha$ is homeomorphic to $[1,\omega^{\omega^\zeta \alpha}]$,
there is a closed neighborhood $M_n$ of $l_n$ in $L_\alpha$ and an ordinal
$\gamma''$, $\omega^\zeta \beta\geq
\gamma''\geq \gamma,$ such that $M_n^{(\gamma'')}=\{l_n\}.$ Moreover, we
may assume that there is a closed subset $M_n^{ \prime}$
 of $L_\beta$ such that 
\[
M_n=\{\frac 12(\sum_{j\in K_\beta }l^{\prime }(j)\delta
_{(k(n),l_n^{^{\prime \prime }},m,j)}+\sum_{j\in K_\beta }l_n^{\prime
}(j)\delta _{(k_1,l_n^{\prime \prime },m,j)}):l^{\prime }\in M_n^{
\prime} \}.
\]
Because $(l_n,G_n)\in P_\gamma (\epsilon ,L_\alpha),$
$(l_n,G_n)\in P_\gamma (\epsilon ,M_n).$ Let
\begin{align*}
M_n^1 &=\{\sum_{j\in K_\beta }l_n^{\prime }(j)\delta
_{(k(n),l_n^{^{\prime \prime }},m,j)}:l_n^{\prime}\in M_n^{\prime}\} \\
\intertext{and}
M_n^2 &=\{\sum_{j\in K_\beta }l_n^{\prime }(j)\delta
_{(k_1,l_n^{^{\prime \prime }},m,j)}:l_n^{\prime}\in M_n^{\prime}\}.
\end{align*}
Then a
transfinite induction argument shows that 
$(l_n^1,G_n^1)\in P_{\gamma }(\epsilon ,M_n^1)$ or $(l_n^2,G_n^2)\in P_{\gamma 
}(\epsilon , M_n^2)$ for infinitely many $n,$ where $%
G_n^{1}=\{j:(k(n),l_n^{\prime \prime },m,j)\in G_n\}$ and $%
G_n^{2}=\{j:(k_1,l_n^{\prime \prime },m,j)\in G_n\}$ (The
argument is essentially the same as the proof of Lemma \ref{coL}.) Because
the mapping $\psi :(k,l,m,j)\rightarrow j$ is a homeomorphism,
$(l_n^i,\psi(G_n^i))\in P_\gamma(\epsilon,M_n^i)$ for $i=1$ or $2.$
Thus we must have $\rho(\epsilon,L_\beta)\geq \gamma.$
In particular, if 
$l\in \Phi (L_1)^{(0)},$ we have that $(l,G)\in P_{\gamma+1} (\epsilon
,L_\alpha )$ only if $P_\gamma (\epsilon ,L_\beta )\neq \emptyset .$

Notice that in the situation above
 if $\phi(k_0)\notin G,$ then $\phi(k_0)\notin \overline{\cup G_n}.$ Thus
for large $n,$ $l_n^2(G_n)=0.$
In order for $l_n(G_n)\geq \epsilon,$ we must have $%
l_n^1(G_n)\geq 2\epsilon ,$ and $(l_n^1,G_n^1)\in P_\gamma(2\epsilon,M_n^1).$
Thus $P_\gamma(2\epsilon ,L_\beta )\neq \emptyset,$ and $(l_n,G_n)\in P_\gamma
(\epsilon ,L_\alpha),$ only if $(l_n',\psi(G_n^1)) \in 
P_{\gamma }(2\epsilon ,L_\beta).$

With these observations we can now estimate $\rho (\epsilon ,\alpha ).$

Suppose that $(l_0,G_0)\in P_\gamma (\epsilon ,L_\alpha )$ for some $\gamma
\geq \max (\rho (2\epsilon ,\beta )+\omega ^\zeta +1,\rho (\epsilon ,\beta
)+2).$ Then there exists $(l_{0,n},G_{0,n})\in P_{\gamma _n}(\epsilon
,L_\alpha ),$ where $\gamma _n=\gamma -1,$ if $\gamma $ is a successor
ordinal or $\gamma _n\uparrow \gamma ,\gamma _n\geq $ $\max (\rho (2\epsilon
,\beta )+\omega ^\zeta ,\rho (\epsilon ,\beta )+1),$ if $\gamma $ is a limit
ordinal, such that $l_{0,n}\stackrel{\text{w}^{*}}{\rightarrow }%
l_0,G_{0,n}\cap G_{0,n^{\prime }}=\emptyset ,$ for $n\neq n^{\prime },$ and $%
\overline{\cup G_{0,n}}\subset G_0.$ If $l_{0,n}\in \Phi (L_1),$ for
infinitely many $n,$ then for at most one $n,$ $\phi (k_1)\in G_{0,n}.$ We
can assume, by discarding that one, that there is no such $n.$ Also in this
case $\epsilon \leq 1/2,$ since $\rho (\epsilon ,1)=0$ for $\epsilon >1/2.$
Because $\gamma \geq \rho (\epsilon ,\beta )+1\geq \omega ^\zeta +1$ and $%
l_{0,n}\notin \Phi (L_1)^{(\omega ^\zeta )},$ there is a w$^{*}$-open
neighborhood $O$ of $l_{0,n},$ such that $P_{\omega ^\zeta }(\epsilon ,O\cap
\Phi (L_1))=\emptyset .$ and $(l_{0,n},G_{0,n})\in P_{\gamma _n}(\epsilon
,O).$ Therefore there exists $((\mu _n,F_n))_{n=1}^\infty \subseteq P_\eta
(\epsilon ,O\setminus \Phi (L_1))$ with $\mu _n\stackrel{\text{w}^{*}}{%
\rightarrow }l_{0,n},$ $\overline{\cup F_i}\subseteq G_{0,n},$ $F_i\cap
F_{i^{\prime }}=\emptyset ,$ for $i\neq i^{\prime }$ and $\eta +\xi \geq
\gamma _n$ for some $\xi <\omega ^\zeta .$ Because $\phi (k_1)\notin
G_{0,n}, $ it follows from the argument above that $P_\eta (2\epsilon
,L_\beta )\neq \emptyset .$ Thus $\eta \leq \rho (2\epsilon ,\beta ).$ This
contradicts the choice of $\gamma \geq \rho (2\epsilon ,\beta )+\omega
^\zeta +1.$ Hence $l_{0,n}\in \Phi (L_1)$ for only finitely many $n.$ This
implies that $P_{\gamma _n}(\epsilon ,L_\beta )\neq \emptyset $ for all but
finitely many $n$ and hence $\gamma _n\leq \rho (\epsilon ,\beta ),$ again a
contradiction of the choice of $\gamma .$ Thus we have that $\rho (\epsilon
,\beta +1)\leq \max (\rho (2\epsilon ,\beta )+\omega ^\zeta ,\rho (\epsilon
,\beta )+1)$ for $\epsilon \leq 1/2.$ If $\epsilon >1/2,l_{0,n}\notin \Phi
(L_1)$ and thus $\gamma _n\leq \rho (\epsilon ,\beta ).$ Hence $\rho
(\epsilon ,\alpha )\leq \rho (\epsilon ,\beta )+1\leq \max (0+0,\rho
(\epsilon ,\beta )+1)=\max (\rho (2\epsilon ,\beta )+\rho (\epsilon ,1),\rho
(\epsilon ,\beta )+1)$ for $\epsilon >1/2.$

In the case that $\alpha $ is a limit ordinal it is easy to see that $\rho
(\epsilon ,\alpha )\leq \sup \{\rho (\epsilon ,\beta ):\beta <\alpha \}$
because we have simply glued the spaces $K_{\alpha _n},n\in \mathbb{N}$ at
their distinguished points to make $K_\alpha .$
\end{proof}

\begin{proposition}
Suppose $\alpha <\omega _1$ and $\alpha =\beta +\eta ,$ where $\eta <\omega
^{\zeta +1}$ and $\beta $ is the smallest ordinal for which there exists
such an $\eta .$ Then $\rho (\epsilon ,\alpha )\leq \beta +\omega ^\zeta
l+\eta $ for $2^{-(l+1)}<\epsilon \leq 2^{-l},$ $l=0,1,2,\ldots .$
\end{proposition}

\begin{proof}
The proof is by induction on $\alpha .$ We have already computed $\rho
(\epsilon ,1)$ and it clearly satisfies the inequality. Suppose that it is
true for all $\alpha ^{\prime }<\alpha .$ Fix $l$ and $\epsilon .$

If $\alpha $ is a limit ordinal and $\alpha _n\uparrow \alpha ,$ let $\alpha
_n=\beta _n+\eta _n$ be the decomposition of $\alpha _n$ with $\beta _n$
minimal and $\eta _n<\omega ^{\zeta +1}.$ If the sequence $(\beta _n)$ is
not eventually constant, then $\alpha =\alpha +0$ is the decomposition of $%
\alpha $ and $\alpha =\lim \alpha _n=\lim \beta _n=\lim \beta _n+\omega
^\zeta l+\eta _n\geq \rho (\epsilon ,\alpha ),$ by Lemma \ref{ind}. If there
is some $n_0$ such that $\beta _n=\beta ^{\prime }$ for all $n\geq n_0,$ $%
\lim \beta _n+\omega ^\zeta l+\eta _n=\beta _0+\omega ^{\zeta }l+\lim \eta
_n.$ If $\lim \eta _n<\omega ^{\zeta +1},$ then $\alpha =\beta _0+\lim \eta
_n$ is the decomposition of $\alpha $ and the inequality holds$;$ if not, $%
\eta =0$ and $\alpha =(\beta _0+\omega ^{\zeta +1})+0$ is the decomposition.
In this case $\omega ^{\zeta }l+\lim \eta _n=\omega ^{\zeta +1}$ and thus $%
\rho (\epsilon ,\alpha )\leq \beta _0+\omega ^{\zeta +1}=\beta .$

If $\alpha =\alpha ^{\prime }+1,$ and $\alpha ^{\prime }=\beta ^{\prime
}+\eta ^{\prime }$ is the decomposition of $\alpha ^{\prime },$ then $\alpha
=\beta ^{\prime }+(\eta ^{\prime }+1)$ is the decomposition of $\alpha .$ By
Lemma \ref{ind}, $\rho (\epsilon ,\alpha )\leq \max (\rho (2\epsilon ,\alpha
^{\prime })+\rho (\epsilon ,1),\rho (\epsilon ,\alpha ^{\prime })+1).$ If $%
\epsilon >1/2,l=0,$ $\rho (2\epsilon ,\alpha ^{\prime })+\rho (\epsilon
,1)=0 $ and $\rho (\epsilon ,\alpha ^{\prime })+1\leq \beta ^{\prime }+\eta
^{\prime }+1,$ as required. If $\epsilon \leq 1/2,$ 
\begin{align*}
\max (\rho (2\epsilon ,\alpha ^{\prime })+\rho (\epsilon ,1),\rho (\epsilon
,\alpha ^{\prime })&+1) \\
&\leq \max (\beta ^{\prime }+\omega ^\zeta
(l-1)+\eta ^{\prime }+\omega ^\zeta ,\beta ^{\prime }+\omega ^\zeta l+\eta
^{\prime }+1) \\
&=\beta ^{\prime }+\omega ^\zeta l+\eta ^{\prime }+1 \\
&=\beta +\omega ^\zeta l+\eta .
\end{align*}
\end{proof}

We now have all of the tools to prove our main result.

\begin{theorem}
$\label{main}$If $1\leq \zeta <\alpha <\zeta \omega <\omega _1,$ then there
is an operator $T$ from $C(\omega ^{\omega ^\alpha })$ onto itself such that
if $Y$ is a subspace of $C(\omega ^{\omega ^\alpha })$ which is isomorphic
to $C(\omega ^{\omega ^\alpha })$ then $T_{|Y}$ is not an isomorphism.
\end{theorem}

\begin{proof}
Let $\gamma $ satisfy $\alpha =\zeta +\gamma .$ For $(K_{\omega ^\gamma
},L_{\omega ^\gamma })$ constructed for the ordinal $\zeta ,$ i.e., $K_1$
is homeomorphic to $[1,\omega ^{\omega ^\zeta }],$ and for $2^{-(l+1)}<\epsilon
\leq 2^{-l},\rho (\epsilon ,\omega ^\gamma )$ $\leq \beta +\omega ^\zeta
l+\eta ,$ where $\omega ^\gamma =\beta +\eta $ and $\eta <\omega ^{\zeta
+1}.$ Therefore $\rho (\epsilon ,\omega ^\gamma )<\beta +\omega ^{\zeta +1}$
for every $\epsilon >0.$ Observe that if $\gamma >1,\omega ^\alpha =\omega
^{\zeta +\gamma }>\max (\omega ^\gamma ,\omega ^{\zeta +1})2\geq $ $\rho
(\epsilon ,\omega ^\gamma )$ and if $\gamma =1,$ $\omega ^{\zeta
}l+\omega ^\gamma \leq \omega ^\zeta (l+1)<\omega ^{\zeta +1}=\omega ^\alpha
.$ Thus by Theorem \ref{AW} and Lemma \ref{coL}, there is no subspace $Y$ of 
$C(\omega ^{\omega ^\alpha })$ which is isomorphic to $C(\omega ^{\omega
^\alpha })$ such that $T_{|Y}$ is an isomorphism.
\end{proof}

The failure of the condition given in Theorem \ref{main} for all $\zeta$
is equivalent
to $\alpha =\omega ^\gamma ,$ for some $\gamma<\omega_1.$
 This is still a long way from Bourgain's condition $\omega \omega ^\alpha
=\omega ^{\omega ^\alpha },$ which guarantees the existence of subspaces
isomorphic to $C(\omega^{\omega^\alpha})$ on which maps of
$C(\omega^{\omega^\alpha})$ onto itself would be isomorphisms.
On the other hand the estimate for the
inductive step given in Lemma \ref{ind} is sometimes generous. In specific cases
it is the case that $\rho (\epsilon ,\beta )=\rho (\epsilon ,\beta +1).
$ Thus it could be that a more careful estimate of $\rho (\epsilon ,\alpha )$
would yield a stronger result.  It also seems likely that
there is room for improvement in Bourgain's estimates.

\end{document}